\documentclass[a4paper,11pt]{article}

\usepackage[T1]{fontenc}
\usepackage[utf8]{inputenc}
\usepackage[english]{babel}
\usepackage{lmodern}
\usepackage{amsmath,amssymb,amsthm,mathrsfs}
\usepackage{xcolor}
\usepackage[margin=28mm]{geometry}
\usepackage{microtype}
\usepackage{varioref}
\usepackage[hidelinks]{hyperref}
\usepackage[nameinlink]{cleveref}
\usepackage{commath}
\usepackage{thm-restate}

\DeclareMathOperator{\RR}{\mathbb{R}}

\DeclareMathOperator{\ZZ}{\mathbb{Z}}

\DeclareMathOperator{\IC}{\mathcal{I}}

\DeclareMathOperator{\rank}{rank}

\DeclareMathOperator{\conv}{conv\!.\!hull}

\DeclareMathOperator{\vertex}{vert}

\DeclareMathOperator{\BUnit}{\mathbf 1}

\newcommand*{\intint}[2][1]{\left\{#1,\, \dots,\, #2\right\}}

\newif\ifhighlightchanges
\highlightchangestrue
\newcommand{\gribHighlightColor}[1]{\ifhighlightchanges\color{#1}\fi}

\newcommand\restr[2]{{
  \left.\kern-\nulldelimiterspace 
  #1 
  \vphantom{\big|} 
  \right|_{#2} 
  }}

\highlightchangesfalse

\newtheorem{theorem}{Theorem}
\newtheorem{lemma}{Lemma}
\newtheorem{proposition}{Proposition}
\newtheorem{corollary}{Corollary}

\theoremstyle{definition}

\theoremstyle{definition}

\theoremstyle{plain}

\title{An upper bound on the number of edges\protect\\
of a Gomory corner polyhedron}
\author{
  A.\,Yu.~Chirkov\thanks{A.\,Yu.~Chirkov and S.\,I.~Veselov are listed first in recognition of their foundational contribution to this work.}, S.\,I.~Veselov,\\ D.\,V.~Gribanov,  D.\,S.~Malyshev,\\   N.\,Yu.~Zolotykh
}
\date{}
\hypersetup{
  pdftitle={An upper bound on the number of edges of a Gomory corner polyhedron},
  pdfauthor={A. Yu. Chirkov, S. I. Veselov, D. V. Gribanov, D. S. Malyshev, N. Yu. Zolotykh},
  pdflang={en}
}

\begin{document}
\raggedbottom
\maketitle

\begin{abstract}
We establish the first nontrivial upper bound on the number of edges of a
Gomory corner polyhedron, resolving a long-standing open question. For an
$n$-dimensional corner polyhedron associated with an integer lattice of
determinant $\delta$, the number of edges is at most
\[
O\left(n6^n \cdot \bigl(1+\log_2(2\delta)\bigr)^{2n-3}\right).
\]
This improves the standard bound obtained by counting pairs of vertices,
reducing the logarithmic exponent from $(2n-2)$ to $(2n-3)$.

\end{abstract}

\section{Introduction}

The \emph{Gomory corner polyhedron} is the convex hull of the set
\[
  G=\bigl\{x\geq 0 : x\in b+L\bigr\},
\]
where $L \subseteq \ZZ^n$ is an $n$-dimensional lattice and $b\in\ZZ^n \setminus L$.

The concept of a \emph{corner polyhedron} was introduced by
R.~E.~Gomory in \cite{GomoryCombinatorialPoly}, where its fundamental
properties were also established. In the same work, Gomory described
a class of so-called \emph{master polyhedra}, showed how their facets can
be used to obtain facets of the corresponding corner polyhedra, and gave
a complete characterization of the facets of a master polyhedron.
This characterization admits an algorithmic implementation, which in turn
provides cutting planes for branch-and-cut algorithms for integer linear
programming \cite{gomory_corner_survey}.
See also Hu's book \cite{HuBook} and the more recent survey by Gribanov
and Malyshev \cite{GomoryGroupSurvey}.

Gomory introduced the concept of an \emph{irreducible point} as a useful
tool for studying the vertices of integer polyhedra: a point $p\in G$
is called irreducible if the inequality $\abs{x}\leq p$ has a unique
solution $x\in L$.

In what follows, however, we use a different definition of
\emph{irreducibility}, which is easily seen to be a special case of
Gomory's definition. We call a point of a set \emph{irreducible} if it
cannot be represented as the midpoint of two other points of that set.
A classical result of V.~Shevchenko \cite{SeparationProperty_origin}
(see also \cite[Theorem~3.13]{BlueBook}) gives the following upper bound
on the number of \emph{irreducible points} with respect to $G$.
\begin{theorem}[Shevchenko~\cite{SeparationProperty_origin}]
\label{thm:irreducible}
Let $X\subseteq G\cap[0,a]^n$. Then the number of irreducible points
in $X$ is at most $(1+\log_2 a)^{n-1}$.
\end{theorem}

Since every vertex of the polyhedron $\conv(G)$ is irreducible and all
vertex coordinates are at most $\delta$\footnote{The latter follows from
the elementary fact that $\delta e_i\in L$ for the $i$th standard basis
vector $e_i$.}, \Cref{thm:irreducible} implies an upper bound on the number
of vertices. As shown in \cite[Section~4]{IntVertEstimates_VesChir_RU2007},
the dependence of this bound on $\delta$ is best possible.
\begin{corollary}\label{cor:G_vertices}
  The number of vertices of the polyhedron $\conv(G)$ is at most
  \begin{equation*}
    (1+\log_2\delta)^{n-1}.
  \end{equation*}
\end{corollary}

We establish the first nontrivial upper bound on the number of
one-dimensional faces of $\conv(G)$ that is polynomial in $\log\delta$
for fixed $n$. Estimating the number of extreme rays (unbounded
one-dimensional faces) is straightforward (see \Cref{th:ray-count}).
However, the only way known to us to estimate the number of edges
(bounded one-dimensional faces) is to count all unordered pairs of
vertices. This gives the bound
\begin{equation}\label{eq:simple_edges_UB}
  (1 + \log \delta)^{2n-2}.
\end{equation}

Improving \eqref{eq:simple_edges_UB} has long been an open problem in
polyhedral combinatorics and integer linear programming. This question
is discussed, for example, by B\'ar\'any, Howe and Lov\'asz
\cite{IntVertLower_BaranyHoweLovasz1992}. Our \emph{main result} is the
first upper bound that improves the exponent in
\eqref{eq:simple_edges_UB}. The proof is given in
\Cref{sec:proof:thm:edges}.
\begin{restatable}{theorem}{ThEdges}\label{thm:edges}
The number of one-dimensional faces of the polyhedron $\conv(G)$ is at most
\[
  2n 6^n \cdot (1+\log_2 2\delta)^{2n-3}.
\]
\end{restatable}


\subsection{Motivation}
\label{sec:motivation}

Our main motivation for studying Gomory corner polyhedra is their role in
obtaining upper bounds on the number of vertices of integer hulls of general
polyhedra. Consider
\[
  P=\{x\in\RR^n:Ax\leq b\},
  \qquad A\in\ZZ^{m\times n},\quad b\in\ZZ^m,
\]
and let $P_I=\conv(P\cap\ZZ^n)$ be its integer hull.
We assume that $P$ is nonempty and pointed; for a nonempty $P$, pointedness
is equivalent to $\rank(A)=n$.
Write
\[
  \Delta=\Delta_n(A),
  \qquad \Delta_{\mathrm{ext}}=\Delta_n\bigl((A\mid b)\bigr),
\]
where $\Delta_n(M)$ denotes the maximum absolute value of the determinants
of all $n\times n$ submatrices of $M$. In particular, the definition of
$\Delta_{\mathrm{ext}}$ uses minors of order $n$, even though the augmented
matrix has $n+1$ columns.

Classical methods for estimating
$\abs{\vertex(P_I)}$ are surveyed by Zolotykh
\cite{Vertices_Survey_Zolotykh} and Veselov and Chirkov
\cite{IntVerticesSurveyPart1,IntVerticesSurveyPart2}.
See also
\cite{IntVertEstimates_VesChir_RU2007,GomoryCombinatorialPoly,IntVertKnapsack_Hayes,IntVert_Morgan1991}
for further background and related results.
An earlier general upper bound, due to Cook, Hartmann, Kannan and McDiarmid
\cite{IntVert_Cook}, is polynomial in the binary encoding length of the
defining system for fixed dimension. Subsequent bounds of Chirkov and
Veselov improve the dependence on the number of inequalities and, in a
further estimate, eliminate dependence on the right-hand side. We therefore
state only these stronger bounds. In terms of the determinant parameters
above, they take the following form (see
\cite{IntVerticesSurveyPart1,IntVerticesSurveyPart2}):
\begin{align*}
  \abs{\vertex(P_I)}
  &=O_n\!\left(m^{\lfloor n/2\rfloor} \cdot 
       (1+\log\Delta_{\mathrm{ext}})^{n-1}\right),
       \\
  \abs{\vertex(P_I)}
  &=O_n\!\left(m^{n} \cdot
       (1+\log\Delta)^{n-1}\right).
\end{align*}
Thus, for fixed dimension, these estimates are polynomial in $m$ and the
logarithm of the relevant determinant parameter. The first bound depends on
$b$ through $\Delta_{\mathrm{ext}}$, whereas the second is independent of
$b$.

B\'ar\'any, Howe and Lov\'asz \cite{IntVertLower_BaranyHoweLovasz1992}
gave a fairly involved procedure for constructing rational simplices
whose integer hulls have
$\Omega_n(\varphi^{n-1})$ vertices, where $\varphi$ is the total binary
encoding length of the defining inequalities. Thus, for fixed dimension,
the exponent $n-1$ in the encoding-length parameter is optimal. Veselov obtained the stronger lower bound
$\Omega_n\bigl((\log\Delta)^{n-1}\bigr)$ for the average number of vertices
in a family of Gomory corner polyhedra with lattice determinant $\Delta$;
see \cite{Vertices_Survey_Zolotykh} and
\cite[Section~4]{IntVertEstimates_VesChir_RU2007}.

A complementary upper bound, useful for small determinants, is due to
Gribanov, Shumilov and Malyshev \cite{Gribanov_SimpleConnection}:
\begin{equation*}
  \abs{\vertex(P_I)}\leq 2\binom{m}{n} \cdot \Delta^{n-1}.
\end{equation*}

The connection with Gomory polyhedra is particularly transparent for
simplicial cones. Choose $n$ linearly independent rows of $A$, indexed by
$I$, and consider $Q_I=\{x\in\RR^n:A_Ix\leq b_I\}$. The invertible affine
map $y=b_I-A_Ix$ identifies its integer hull with
\[
  \conv\bigl(\{y\in\RR_{\geq0}^n \colon y\in b_I+A_I\ZZ^n\}\bigr).
\]
This is a Gomory corner polyhedron with lattice determinant
$\abs{\det A_I}\leq\Delta$, except in the trivial case
$b_I\in A_I\ZZ^n$.
Bounds for these basic corner polyhedra, together with suitable covering
and reduction arguments, are ingredients in estimates for general integer
hulls.

\begin{quote}
  \it This approach motivates an analogous programme for the number of edges
and higher-dimensional faces: combine bounds for Gomory corner polyhedra
with suitable reductions to obtain bounds for general integer hulls.
Before developing such general bounds, we therefore seek the strongest
possible estimates for Gomory corner polyhedra themselves.
The improvement of the edge bound in \Cref{thm:edges} is a step in this
direction.
\end{quote}

\subsection{Notation and definitions}

We denote the determinant of the lattice $L$ by $\delta$.
A point $p$ of a set is called \emph{irreducible} if it cannot be represented
as the midpoint of two other points of that set.

A vector $\pi\in\RR_{\geq 0}^n$ is called a \emph{supporting vector}
for $G$ if $\pi^\top x\geq 1$ for every $x\in G$ and the set
\[
  \Pi=\{x\in G : \pi^\top x=1\}
\]
is non-empty. The set $\Pi$ is called a \emph{face} of the polyhedron
$\conv G$.
A \emph{vertex} is a face of affine dimension zero; an \emph{edge} and
an \emph{extreme ray} are faces of affine dimension one. The difference
is that an edge is determined by two vertices, whereas an extreme ray is
determined by a vertex and a direction vector. An edge with vertices
$p$ and $q$ is denoted by $\{p,q\}$.

We set
\begin{gather*}
  M(a) = [0,1+\log_2 a]^n,\\
  B(t) =\bigl\{x\in\ZZ^n : 2^{t_i}-1\leq x_i\leq 2^{t_i+1}-2,
                \ i \in \intint n\bigr\}.
\end{gather*}
Let $V$ and $E$ be the sets of vertices and edges of $\conv(G)$, respectively.
The unique box containing a vertex $v$ is denoted by $B(t_v)$. For an edge
$e=\{p,q\}$, the unique box containing the sum $p+q$ of its endpoints is
denoted by $B(t_e)$. We set
\[
  E(t)=\{e\in E : t_e=t\}.
\]

\section{An upper bound on the number of extreme rays}

\begin{lemma}
\label{lem:ray-coordinates}
Let $r$ be an extreme ray emanating from a vertex $v$, and let $v_r$ be
the point with integer coordinates on $r$ closest to $v$.
The coordinates of $v_r$ do not exceed $2\delta$.
\end{lemma}

\begin{proof}
Let $\pi$ be a supporting vector for the ray $r$, and put $h=v_r-v$.
By the choice of $v_r$, $h$ is an integer direction vector of $r$ of minimum length.
Clearly, $\pi \geq 0$, $h\geq 0$, and $\pi^\top h=0$.
Consequently, $\pi_i h_i=0$ for every $i$.
Suppose that $h_i>\delta$ for some $i$. Consider the vector
$h'=h-\delta e_i$, where $e_i$ is a standard basis vector in $\RR^n$.
Since $h'\geq 0$, $\pi^\top\cdot h'=0$, and $h'\in\Lambda$, the vector $h'$
is a direction vector of the ray $r$. However, its length is smaller
than that of $h$, a contradiction.
\end{proof}

\begin{proposition}
\label{th:ray-count}
The number of extreme rays does not exceed $(1+\log_2 2\delta)^{n-1}$.
\end{proposition}

\begin{proof}
Let $V'$ be the set of integer points, each lying on one of the rays and closest
to the vertex of that ray. Clearly, every point of $V\cup V'$ is
irreducible.
\end{proof}

\section{An upper bound on the number of edges (proof of \texorpdfstring{\Cref{thm:edges}}{Theorem \getrefnumber{thm:edges}})}\label{sec:proof:thm:edges}

\noindent The proof of \Cref{thm:edges} requires several auxiliary lemmas, which we establish first. At the end of this section, we will derive the theorem directly from these lemmas.

\begin{lemma}
\label{lem:vertex-edge}
Let $v\in V$, $e=\{p,q\}\in E$, and $t_v<t_e$.
Then $v\in\{p,q\}$.
\end{lemma}

\begin{proof}
Put $v'=p+q-v$. Clearly, $v'\in b+L$.
The inequalities
\[
  v\leq 2^{t_v+1}-2 < 2^{t_e}-1\leq p+q
\]
imply that $v'\geq 0$. Consequently, $v'\in G$.
Let $\pi$ be a supporting vector for $e$. Since $\pi^\top v\geq 1$,
$\pi^\top v'\geq 1$, and
\[
  \pi^\top v+\pi^\top v'=\pi^\top p+\pi^\top q=2,
\]
we have $\pi^\top v=1$. Thus, the vertex $v$ lies in a supporting hyperplane
of $G$ containing only the vertices $p$ and $q$, i.e., $v\in\{p,q\}$.
\end{proof}

\begin{corollary}
\label{cor:common-vertex}
Suppose that $t_p<t$ for a vertex $p$ of one of the edges in $E(t)$.
Then $p$ is a vertex of every edge in $E(t)$.
\end{corollary}

\begin{lemma}
\label{lem:edges-in-box}
For every $t\in M(2\delta)$,
\[
  \abs{E(t)} \leq 2 \cdot 6^n \cdot (1+\log_2\delta)^{n-2}.
\]
\end{lemma}

\begin{proof}
    Let $\{p,q\}$ be an arbitrary edge in $E(t)$. By definition, $p+q\in B(t)$, $p\in B(t_p)$, and $q\in B(t_q)$. The lower bound defining $B(t)$ and the upper bounds defining the endpoint boxes give, for every $i$,
    \begin{gather*}
      2^{t_i}-1\leq p_i+q_i
      \leq (2^{(t_p)_i+1}-2)+(2^{(t_q)_i+1}-2),\quad\text{and}\\
      2^{t_i-1}+2\leq 2^{(t_p)_i}+2^{(t_q)_i}.
    \end{gather*}
    For $t_i \geq 2$, this implies
    \begin{equation}\label{eq:tpq_max_1}
      \max\left\{(t_p)_i, (t_q)_i\right\}\geq t_i-1.
    \end{equation}
    For $t_i\in\{0,1\}$, it follows from the nonnegativity of the box indices. 

    Fix the sets
    \[
      \IC_p=\left\{i\colon (t_p)_i\leq t_i-2\right\},
      \qquad
      \IC_q=\left\{i\colon (t_q)_i \leq t_i-2\right\}.
    \]
    
    We begin by considering the case in which $\IC_p \neq \emptyset$ and $\IC_q \neq \emptyset$. Note that for $i\notin \IC_p$ the coordinate $(t_p)_i$ can take only two values: $t_i-1$ and $t_i$. Consequently, by \Cref{thm:irreducible}, the vertex $p$ can be chosen in at most $2^{n-\abs{\IC_p}} \cdot (1+\log_2\delta)^{\abs{\IC_p}-1}$ ways. The number of choices for $q$ is estimated similarly. By \eqref{eq:tpq_max_1}, we have $\IC_p\cap \IC_q=\varnothing$. Thus,
    \begin{equation*}
      \abs{E(t)} \leq 2^{2n-\abs{\IC_p}-\abs{\IC_q}} \cdot 
              (1+\log_2\delta)^{\abs{\IC_p}+\abs{\IC_q}-2} \leq
      2^n \cdot (1+\log_2\delta)^{n-2}.
    \end{equation*}

    Now consider the case when $\IC_p = \emptyset$. If $\abs{\IC_q} \leq n-1$, then by the same reasoning we get 
    \begin{equation*}
      \abs{E(t)} \leq 2^{2n-\abs{\IC_q}} \cdot (1+\log_2\delta)^{\abs{\IC_q}-1} \leq
      2^{n+1} \cdot (1+\log_2\delta)^{n-2}.
    \end{equation*}

    Therefore, it remains to consider the case $\IC_p = \emptyset$ and $\abs{\IC_q} = n$. The condition $\abs{\IC_q} = n$ implies $t_q < t$. By \Cref{cor:common-vertex}, the vertex $q$ is a vertex of every edge in $E(t)$. 
    If $t_p < t$, then, again by \Cref{cor:common-vertex}, the edge $\{p,q\}$ is the unique edge in $E(t)$. Otherwise, there exists a component $i$ such that $(t_p)_i = t_i$. Thus, the value of $\abs{E(t)}$ is bounded by the number of endpoints $p$ with one fixed coordinate. By \Cref{thm:irreducible}, we get
    \begin{equation*}
      \abs{E(t)} \leq n \cdot (1 + \log_2 \delta)^{n-2}. 
    \end{equation*}

    Therefore, for each fixed configuration of the sets $\IC_p$ and $\IC_q$, we have proven that
    \begin{equation*}
      \abs{E(t)} \leq 2^{n+1} \cdot (1+\log_2\delta)^{n-2}.
    \end{equation*}
    The number of pairs $(\IC_p, \IC_q)$ with $\IC_p \cap \IC_q=\varnothing$ is $3^n$ since each element can belong to $\IC_p$, $\IC_q$ or neither. This gives the desired bound on $\abs{E(t)}$ and completes the proof.
\end{proof}

\begin{lemma}
\label{lem:comparable-boxes}
If $t'<t''$, then at least one of the sets $E(t')$ and $E(t'')$ is empty.
\end{lemma}

\begin{proof}
Suppose, to the contrary, that
$e'=\{p',q'\}\in E(t')$ and $e''=\{p'',q''\}\in E(t'')$.
Let $\pi$ be a supporting vector for $e''$. Since
\[
  2\leq\pi^\top p'+\pi^\top q'\leq\pi^\top p''+\pi^\top q''=2,
\]
we have $\pi^\top p'=\pi^\top q'=1$. Consequently, $e'=e''$ and $t'=t''$,
contrary to the assumption of the lemma.
\end{proof}

\textbf{Let us return to the proof of \Cref{thm:edges}.} For convenience, we recall its statement.
\ThEdges*

\begin{proof}
  Consider the following equivalence relation on the set $M(2\delta)$: $t' \sim t''$ if and only if there exists $\alpha \in \ZZ$ such that
  \[
    t' - t'' = \alpha \cdot \BUnit.
  \]
  Clearly, the number of equivalence classes does not exceed the number of vectors having at least one zero coordinate, i.e., it is at most
  $n \cdot (1+\log_2 2\delta)^{n-1}$. By \Cref{lem:comparable-boxes}, each equivalence class contains at most one $t$ such that $E(t) \not= \emptyset$. The cardinality of each such $E(t)$ is bounded by \Cref{lem:edges-in-box}. Since every edge lies in some $E(t)$, this completes the proof.
\end{proof}

\section*{Acknowledgements}

This result was obtained through joint work at the Nizhny Novgorod Seminar
on Discrete Mathematics. The seminar was founded in 1963 by Yuri
Vasilyevich Glebsky and Alexander Alexandrovich Markov.


\section*{Statements and Declarations}

\subsection*{Competing interests}
The authors declare that they have no competing interests.

\subsection*{Data availability}
The manuscript has no associated data.

\subsection*{Use of artificial intelligence}
The proofs were developed without the use of AI theorem-proving tools; AI tools were used only to assist with surveying the literature and checking for language errors.

\bibliographystyle{unsrt}
\bibliography{grib_biblio}

\end{document}